\voffset-1cm
\input amstex
\magnification=1200
\documentstyle{amsppt}
\NoRunningHeads
\NoBlackBoxes
\topmatter
\title The laced interactive games\linebreak
and\linebreak their {\sl a posteriori\/} analysis
\endtitle
\author Denis V. Juriev
\endauthor
\affil ul.Miklukho-Maklaya 20-180, Moscow 117437 Russia\linebreak
E-mail: denis\@juriev.msk.ru
\endaffil
\date math.OC/9901043\qquad\qquad\qquad\qquad\qquad January 10, 1999\enddate
\endtopmatter
\document

The mathematical formalism of interactive games, which extends one
of ordinary differential games [1] and is based on the concept of an
interactive control, was proposed by the author [2] to take into account 
the complex composition of controls of a real human person, which are often 
complicated couplings of his/her cognitive and known controls with the 
unknown subconscious behavioral reactions of the human organism of this 
person on them.

{\sl A posteriori\/} analysis of the interactive games is of an interest 
as from the internal point of view as providing a way to unravel the 
structure of the so-called {\it intention fields}, one of the main concept 
of the field-theoretic description of the interactive games [2], which 
interactions are explicated by the game. To my mind one of the important 
examples of an intention field is provided by the so-called {\it ``chi''\/} 
of the Chineese tradition, intensively explored in Chi Kung or Feng Shui
(see e.g. the book [3], where the exposition of the main basic ideas of Chi 
Kung is presented in the most systematic form). On the other hand an 
investigation of the structure of intention fields is very important for the 
elaboration of systems of the accelerated nonverbal cognitive computer and 
telecommunications based on the using of droems and their dynamical 
reconstruction [4]. In this case intention fields realise a ``virtual 
dynamical memory'', where the droems are really situated. And, at last, the 
unification of researches in the classical areas of Chi Kung (combined with
Feng Shui, which is specially adapted to the design of the interactive virtual 
realities -- the ``virtual Feng Shui''), nonverbal cognitive communications 
and tele{\ae}sthesy (especially, in presence of external stimuli) with computer 
technologies is very perspective for the development of modern information 
technologies. In all cases mathematical {\sl a posteriori\/} analysis of 
interactive games is truly the most effective and simultaneously inexpensive 
way to understand theoretically the collected experimental data and to make 
them applicable in practice. 

This article is devoted to some aspects of {\sl a posteriori\/} analysis
of certain class of interactive games, namely, the laced interactive 
games, where it may be performed rather effectively. Such analysis allows
to give some approximations of the interactive game by ordinary 
differential games in real time; the obtained series of the approximating
ordinary differential games may be used for the formulation of predictions
for processes in the considered interactive game.

\head\S1. The laced interactive games\endhead

In this paragraph we shall briefly introduce some general concepts of the 
theory of interactive games.

\subhead 1.1. The differential interactive games\endsubhead

\definition{Definition 1 [2]} {\it An interactive system\/} (with $n$
{\it interactive controls\/}) is a control system with $n$ independent controls 
coupled with unknown or incompletely known feedbacks (the feedbacks, which are
called the {\it behavioral reactions}, as well as their couplings with 
controls are of a so complicated nature that their can not be described 
completely). {\it An interactive game\/} is a game with interactive controls 
of each player.
\enddefinition

Below we shall consider only deterministic and differential interactive
systems. For symplicity we suppose that $n=2$. In this case the general
interactive system may be written in the form:
$$\dot\varphi=\Phi(\varphi,u_1,u_2),\tag1$$
where $\varphi$ characterizes the state of the system and $u_i$ are
the interactive controls:
$$u_i(t)=u_i(u_i^\circ(t),\left.[\varphi(\tau)]\right|_{\tau\le t}),$$
i.e. the independent controls $u_i^\circ(t)$ coupled with the feedbacks on
$\left.[\varphi(\tau)]\right|_{\tau\le t}$. One may suppose that the
feedbacks are integrodifferential on $t$ in general, but below we shall 
consider only differential dependence. It means that 
$$u_i(t)=u_i(u_i^\circ(t),\varphi(t),\dot\varphi(t),\ddot\varphi(t),\ldots,
\varphi^{(n)}(t)).$$
A reduction of the general integrodifferential case to the differential one
via an introduction of the intention fields was considered in [2]. 

\subhead 1.2. The phase lacing integrals\endsubhead

\definition{Definition 2} Let us consider a differential interactive game
with two players, $\varphi$ is the state of the system, $u_i$ are
the interactive controls whereas $u^\circ_i(t)$ are the pure controls 
($i=1,2$). The magnitude $K=K(\vec u(t),\vec u^\circ(t),\varphi(t),
\dot\varphi(t))$ is called the {\it phase lacing integral\/} iff it 
is constant in time for all possible parties of the game. Here
$\vec u(t)=(u_1(t),u_2(t))$, $\vec u^\circ(t)=(u_1^\circ(t),u_2^\circ(t))$.
\enddefinition

The definition of the phase lacing integral may be generalized on games
with an arbitrary number of players.

Let us consider some special classes of the phase lacing integrals:
\roster
\item"--" {\it Configuration lacing integrals:\/} $K=K(\vec u(t),\vec 
u^\circ(t),\varphi(t))$. Such phase lacing integrals will be called
simply {\it lacing integrals}.
\item"--" {\it Dynamical lacing integrals:\/} $K=K(\vec u(t),\varphi(t),
\dot\varphi(t))$.
\endroster

\remark{Remark 1} The evolution equations of the interactive game supply
us by the dynamical lacing integrals $K_i=\dot\varphi_i(t)-\Phi_i(\varphi,
\vec u)$, where indices $i$ denote the components of the magnitudes in
any coordinate system.
\endremark

\remark{Remark 2} The phase lacing integrals are closed under natural
operations (summation, multiplication on any real number or on other
phase lacing integrals, functional transformations, etc). The configuration
and dynamical lacing integrals possess the same property.
\endremark

The phase lacing integrals may be considered for any differential interactive 
system. They are analogs of ordinary integrals for dynamical systems.

\subhead 1.3. The laced interactive games\endsubhead

\definition{Definition 3} Let us consider a differential interactive game
with two players, $\varphi$ is the state of the system, $u_i$ are
the interactive controls whereas $u^\circ_i(t)$ are the pure controls 
($i=1,2$), each $u_i$ and $u_i^\circ$ has $n$ degrees of freedom.
The game is called the {\it laced interactive game\/} iff it admits
$2n$ functionally independent over $\vec u$ phase lacing integrals
$K_\alpha(\vec u,\vec u^\circ,\varphi,\dot\varphi)$ ($\alpha=1,2,\ldots,2n$). 
\enddefinition

The functional independence of $K_\alpha$ over $\vec u$ means that
for any fixed values of $\vec u^\circ$, $\varphi$ and $\dot\varphi$ 
the magnitudes $K_\alpha$ are (locally) functionally independent as 
functions of $\vec u$. 

Here and below we shall suppose that the phase lacing integrals depend
smoothly on their arguments. In this situation the (local) functional
independence is equivalent to the inequality of the Jacobian of the
mapping $\vec u\mapsto(K_1,\ldots,K_{2n})$ to zero.

Note that the evolution equations provide us by $m$ dynamical lacing
integrals (see remark 1), where $m$ is the number of degrees of freedom
of the system, i.e. the number of coordinates, which describe the state
$\varphi$. If all such dynamical lacing integrals are functionally
independent we should look for at least $2n-m$ other phase lacing integrals.
Often we are able to choose them from the configuration lacing integrals.

\head\S2. A posteriori analysis of the laced interactive games\endhead

\subhead 2.1. {\sl A posteriori\/} determination of feedbacks\endsubhead

{\sl A posteriori\/} determination of feedbacks means the expression of
$\vec u$ via $\vec u^\circ$, $\varphi$ and $\dot\varphi$ using the
phase lacing integrals $K_{\alpha}$. These determination presupposes the 
knowledge of such {\sl a posteriori\/} data as $\dot\varphi$.

\proclaim{Theorem 1} Let us consider a laced interactive game, $\varphi$ is 
the state of the system, $u_i$ are the interactive controls whereas 
$u^\circ_i(t)$ are the pure controls of two players ($i=1,2$), 
each $u_i$ and $u_i^\circ$ has $n$ degrees of freedom, $K_1,\ldots
K_{2n}$ are the phase lacing integrals. In this case the interactive
controls $u_i$ may be expressed (locally) via pure controls $u_i^\circ$, 
the state $\varphi$, its time derivative $\dot\varphi$ and the phase
lacing integrals $K_\alpha$ as known constants.
\endproclaim

\demo{Proof} One should use the fact that the phase lacing integrals are 
known constants and that the Jacobian of the mapping $\vec u\mapsto(K_1,
\ldots,K_{2n})$ is not equal to zero.
\enddemo

\remark{Remark 3} Note that each interactive control $u_1$ and $u_2$ is
expressed via both $u_1^\circ$ and $u_2^\circ$.
\endremark

\subhead 2.2. Virtual {\sl a posteriori\/} decomposition of a collective 
control\endsubhead

The virtual {\sl a posteriori\/} decomposition of a collective control
means such simultaneous transformation of both interactive and pure
controls of two players to the controls of the virtual players that 
{\sl a posteriori\/} determined feedback of each such player does not contain
a dependence on the controls of the other player.

\proclaim{Theorem 2} Let us consider a laced interactive game, $\varphi$ is
the state of the system, $u_i$ are the interactive controls whereas 
$u^\circ_i(t)$ are the pure controls of two players ($i=1,2$), 
each $u_i$ and $u_i^\circ$ has $n$ degrees of freedom, $K_1,\ldots
K_{2n}$ are the phase lacing integrals. Let the interactive controls $u_i$ 
be expressed (locally) via pure controls $u_i^\circ$, the state $\varphi$, 
its time derivative $\dot\varphi$ and the phase lacing integrals $K_\alpha$ 
as known constants:
$$u_i=u_i(\vec u^\circ,\varphi,\dot\varphi; K_\alpha).$$
Suppose that $\vec u^\circ=0$ implies that $\vec u=0$ (i.e. $0$ is the
stationary point of the mapping $\vec u^\circ\mapsto\vec u$) and also the
Jacobi matrix of the mapping $\vec u^\circ\mapsto\vec u$ is nondegenerate
and diagonalizable at certain neighbourhood of zero then there exist the 
functions $\zeta(\vec x;\varphi,\dot\varphi,K_\alpha)$ locally at some 
neighbourhood of zero so that if we transform both interactive and pure 
controls $\vec u$ and $\vec u^\circ$ according to them:
$$\aligned\vec w=&\zeta(\vec u;\varphi,\dot\varphi,K_\alpha)\\
\vec w^\circ=&\zeta(\vec u^\circ;\varphi,\dot\varphi,K_\alpha)\endaligned$$
then in new variables $\vec w$ and $\vec w^\circ$ {\sl a posteriori\/} 
determined feedbacks would contain the dependences only between the related
components of the controls $\vec w$ and $\vec w^\circ$, i.e.
$$\aligned w_1=&w_1(w_1^\circ,\varphi,\dot\varphi;K_\alpha),\\ 
 w_2=&w_2(w_2^\circ,\varphi,\dot\varphi;K_\alpha).\endaligned$$
\endproclaim

\demo{Proof} The function $\zeta$ is constructed explicitely starting from
the point $\vec u^\circ=0$ to provide the claimed conditions on the
mapping $\vec w^\circ\mapsto\vec w$.
\enddemo

\remark{Remark 4} The controls $w_i$ and $w_i^\circ$ ($i=1,2$) may be
interpreted as interactive and pure controls of two virtual players.
The transformation of the real players to the virtual players depends
on the state $\varphi$ of the game and {\sl a posteriori\/} determined
its time derivative $\dot\varphi$. 
\endremark

\remark{Remark 5 (some psychological interpretations)} The virtual {\sl a 
posteriori\/} decomposition of a collective control may be regarded as an 
identification of two virtual players in the area of collective subconscious 
behavioral reactions. {\sl A posteriori\/} analysis of the laced interactive 
games allows to make such procedure on the mathematical level of precision.
It is remarkable that such decomposition depends on the state $\varphi$
and on its time derivative $\dot\varphi$ (hence, on the intensity of the
real controls). 
\endremark

\head\S3. The retarded control approximation\endhead

{\sl A posteriori\/} analysis of the laced interactive games allows to give 
some approximations of the interactive game by ordinary differential games 
in real time; the obtained series of the approximating ordinary differential 
games may be used for the formulation of predictions for processes in the 
considered interactive game.

\subhead 3.1. The frozen feedback approximation\endsubhead

The simplest approximation of a laced interactive game by the ordinary
differential games is the frozen feedback approximation. Let us consider
the fixed moment $t_0$ of time and the observed values of $\varphi(t_0)$ 
and $\dot\varphi(t_0)$. Then the knowledge of the phase lacing integrals
$K_\alpha$ allows us to express the interactive controls $\vec u(t_0)$ via 
the pure controls $\vec u^\circ(t_0)$:
$$\vec u(t_0)=\vec u(\vec u^\circ(t_0),\varphi(t_0),\dot\varphi(t_0);
K_\alpha).$$
This {\sl a posteriori\/} determined feedback may be frozen and substituted 
into the evolution equations:
$$\dot\varphi(t)=\Phi(\varphi(t),\vec u(\vec u^\circ(t),\varphi(t_0),
\dot\varphi(t_0);K_\alpha)).$$
Thus we obtain an ordinary differential game, which is just the {\it
frozen feedback approximation\/} of the initial laced interactive game.

\subhead 3.2. The retarded control approximation\endsubhead

The retarded control approximation of a laced interactive game is constructed 
in the following manner. Let us consider an arbitrary $\Delta t>0$. For
any moment of time $t$ the knowledge of the phase lacing integrals $K_\alpha$
allows us to perform {\it a posteriori\/} determination of feedbacks, i.e.
to express the interactive controls $\vec u(t)$ via the pure controls
$\vec u^\circ(t)$:
$$\vec u(t)=\vec u(\vec u^\circ(t),\varphi(t),\dot\varphi(t);K_\alpha).$$
Let us now substitute the state $\varphi(t)$ and its time derivative
$\dot\varphi(t)$ by their retarded (delayed) values $\varphi(t-\Delta t)$ and
$\dot\varphi(t-\Delta t)$. Such approximated feedback may be substituted into
the evolution equations:
$$\dot\varphi(t)=\Phi(\varphi(t),\vec u(\vec u^\circ(t),\varphi(t-\Delta t),
\dot\varphi(t-\Delta t);K_\alpha)).$$
Thus we obtain an ordinary differential game {\it (with the retarded, delayed 
arguments)}, which is just the {\it retarded control approximation\/} of the 
initial laced interactive game.

Sometimes the frozen feedback approximation and the retarded control 
approximation coincide but it is not so in general.

\remark{Remark 6} The retarded control approximations form a series
of approximation defined by the initial time $t_0$, such that the states
$\varphi(t)$ for the considered laced interactive game at $t<t_0$ are
the initial data the approximation. One may describe more subtle effects
considering the mutual correlations of the retarded control approximations
with variable $t_0$. 
\endremark

I suspect that the ideological intuition of the nonlinear geometric algebra 
[5] (see also [6]) may be essential for the analysis of series of the retarded 
control approximation, for instance a verification of any algebraic 
correlations between the different retarded control approximations may be 
important as during analysis of concrete games as for the axiomatic separation 
of the interesting classes of the laced interactive games or for a 
classification of various types of interactivity.

Note that the procedure of {\it virtualization\/} may be considered not only 
for the interactive systems [7] but also for the interactive games. 

\remark{Remark 7} The approximating series of games may be used for a
formulation of predictions for processes in the initial interactive game.
Such formulation may be analytic, however, it is difficult to perceive and 
to interpret the obtained results in real time. Thus, it is rather reasonable
to use some {\it visual representation\/} for the series of the approximating
games. Thus, we are constructing an enlargement of the interactive game,
in which the players {\it interactively\/} observe the visual predictions for 
a game in real time. Certainly, such enlargement may strongly transform the
structure of interactivity of the game (i.e. to change the feedbacks entered 
into the interactive controls of players). Note that the aim of the 
virtualization is to restore the past of the interactive processes whereas the 
goal of the proposed enlargement is to correct their future. In some sense
they are two complementary faces of the unique general procedure.
\endremark

\head\S4 Conclusions\endhead

Thus, an important class of differential interactive games, namely, one
of the laced interactive games, was considered. {\sl A posteriori\/} analysis
of such games (including the virtual {\sl a posteriori\/} decomposition
of a collective control) was discussed. Approximations of the laced 
interactive games by the ordinary differential games, the frozen feedback 
approximation and the retarded control approximation, were constructed.

\Refs
\roster
\item"[1]" Isaaks R., Differential games. A mathematical theory with
applications to warfare and pursuit, control and optimization. Wiley,
New York, 1965;\newline
Owen G., Game theory, Saunders, Philadelphia, 1968.
\item"[2]" Juriev D., Interactive games and representation theory. I,II.
E-prints: math.FA/9803020, math.RT/9808098.
\item"[3]" Juriev D., Droems: experimental mathematics, informatics and
infinite dimensional geometry. Report RCMPI-96/05$^+$ [e-version:
cs.HC/9809119].
\item"[4]" Yang Jwing-Ming, The root of Chineese Chi Kung. The secrets
of Chi Kung training. YMAA, Boston, 1989.
\item"[5]" Sabinin L.V., On the nonlinear geometric algebra [in Russian]. In
``Webs and quasigroups''. Kalinin [Tver'], 1988, pp.32-37;\newline
Sabinin L.V., Differential geometry and quasigroups [in Russian]. 
Trans.Inst.Math. Si\-be\-rian Branch Soviet Acad.Sci., Novosibirsk, 1984, v.14, 
pp.208-221.
\item"[6]" Mikheev P.O., Sabinin L.V., Quasigroups and differential geometry. 
In ``Quasigroups and loops. Theory and applications''. Berlin, Heldermann 
Verlag, 1990, P.357-430.
\item"[7]" Juriev D., On the description of physical interactive information
systems [in Russian]. Report RCMPI-96/05 (1996) [e-version: mp\_arc/96-459].
\endroster
\endRefs
\enddocument